\documentclass{alf-j-l}
\usepackage{verbatim} 
\usepackage{syntonly}
\usepackage{url}  

\ifx\preambleloaded\relax
   \else\let\preambleloaded\relax\fi
\def\Z{{\mathbb Z}} \def\Q{{\mathbb Q}} 

\def\F{{\mathbb F}} 

\def\={\equiv}

\def\\{\cr} 
\def\lt{<}
\def\gt{>}
\def\conj{\overline}
\def\ov{\overline}

\def\acclap#1{\raise\hgtsig\hbox to0pt{$#1$\hss}}
\newdimen\hgtsig
\setbox0=\hbox{$\displaystyle{\sum}$}
\hgtsig=\ht0\relax
\advance\hgtsig by -1.75ex

\def\id#1{\langle #1\rangle} 

\newbox\boxW\newdimen\dimW
\def\heighten#1{%
\setbox\boxW\hbox{$\displaystyle #1$}
\dimW=1.04\ht\boxW\advance\dimW by 1.00pt
\vbox to \dimW{}}

\DeclareMathAlphabet{\Bi}{OT1}{cmm}{b}{it}  

\providecommand{\bysame}{\makebox[3em]{\hrulefill}\thinspace}
\def\con#1=#2(#3){#1\equiv#2\pod{#3}}
                            {\end{enumerate}}
                            {\end{enumerate}}

\def\CF#1{{\def\\{\mathrel{,}}\def\;{\mathrel{;}}
\def\dots{\ldots}\def\dotss{\ldots\ldots}[\,#1\,]}}

\def\cfraci#1#2{#1_0+{
\let\ds\displaystyle
\def\vl{\ds1\vrule width0pt depth.5ex height2ex\over}
\vl{\ds#1_1 + {\vl{\ds#1_2 +{\vl\hskip.5em\ddots}}}}}}

\def\divides{{\mathchoice{\mathrel{\bigm|}}{\mathrel{\bigm|}}{\mathrel{|}}%
{\mathrel{|}}}}
\def\Div{\divides}
\def\notdivides{\mathrel{\kern-3pt\not\!\kern3.5pt\bigm|}}

\newbox\boxA
\newbox\boxB
\newdimen\dimA
\newdimen\dimB
\newdimen\dimC

\def\house#1{
\setbox\boxA\hbox{$\displaystyle #1$}
\dimA=1.04\wd\boxA\advance\dimA by 2pt
\dimB=1.04\ht\boxA\advance\dimB by 2pt
\dimC=0.05\wd\boxA
\hskip\dimC\hskip1pt\hbox to \dimA
{\vrule\vbox to \dimB{\hsize=\dimA
\hrule\vfill \centerline{\box\boxA}}\vrule}
\hskip 0.8pt\hskip1pt\hskip\dimC}

\def\poly {polynomial}

\def\cf {continued fraction}
\def\pq {partial quotient}
\def\cq {complete quotient}

\def\ex{expansion}

\def\cfe{continued fraction expansion}

\newcount\hours
\newcount\minutes
\def \SetTime{\hours=\time
\global\divide\hours by 60
\minutes=\hours
\multiply\minutes by 60
\advance\minutes by-\time
\global\multiply\minutes by-1 }
\SetTime
\def \now{\number\hours:\ifnum\minutes<10 0\fi\number\minutes}
\def \Now{\today\ $[$\now$]$}

\newif\ifMacTextures
\MacTexturestrue

\ifMacTextures

\gdef\EPSF#1by#2(#3){%
\vbox to #2{\hrule width #1 height 0pt depth 0pt%
\vfill\special{illustration #3}}}%

\gdef\scaledEPSF#1by#2(#3 scaled #4){{%
\dimen0=#1 \dimen1=#2%
\divide\dimen0 by 1000 \multiply\dimen0 by #4%
\divide\dimen1 by 1000 \multiply\dimen1 by #4%
\EPSF \dimen0 by \dimen1 (#3 scaled #4)}}%
\else
\input epsf

\gdef\EPSF#1by#2(#3){%
\vbox to #2{\hrule width #1 height 0pt depth 0pt\vfill \epsfbox{#3}}}%

\gdef\EPSF#1by#2(#3){\epsfbox{#3}}%

\gdef\scaledEPSF#1by#2(#3 scaled #4){{%
\dimen0=#1 \dimen1=#2%
\divide\dimen0 by 1000 \multiply\dimen0 by #4%
\divide\dimen1 by 1000 \multiply\dimen1 by #4%
\epsfxsize=\dimen0\epsfbox{#3}}}%
\fi
\MacTexturesfalse

\def\ds{\displaystyle}

\theoremstyle{plain}  
\newtheorem{theorem}{Theorem}

\newtheorem{proposition}[theorem]{Proposition} 

\newtheorem*{theorem*}{Theorem}

\theoremstyle{definition}

\newtheorem*{Remark*}{Remark}
\newtheorem*{Remarks*}{Remarks}

\newtheorem*{example*}{Example} 
\newtheorem*{guess*}{Guess}

\theoremstyle{remark}

\newtheorem*{remark*}{Remark} 
\newtheorem*{remarks*}{Remarks}

\newtheoremstyle{aside}
   {6pt}
   {6pt}
   {\footnotesize}
   {}
   {\scshape}
   {:}
   {.5em}
   {}

\theoremstyle{aside}
\newtheorem*{aside}{Aside}
\newtheorem*{aside*}{Concluding Aside}

\PII{Centre for Number Theory Research, 1 Bimbil Place, Killara,
Australia 2071}%

\copyrightinfo{\number\year}{Alfred J van der Poorten}

\begin{document}

\def\currentvolume{169}
\def\currentissue{Draft of\ }
\def\paperdate{\today}
\def\ISSN{}
\pagespan{69}{\pageref{page:lastpage}}

\title {Elliptic
Curves\\ and 
Continued Fractions} 

\author{Alfred J. van der Poorten}
\address{Centre for Number Theory Research, 1 Bimbil Place, Killara, Sydney,
NSW 2071, Australia}
\email{alf@math.mq.edu.au (Alf van der Poorten)}

\thanks{The author was supported in part by a grant from the Australian
Research Council.}

\subjclass[2000]{Primary: 11A55, 11G05; Secondary: 14H05, 14H52}

\date{\Now.}

\keywords{continued fraction expansion, function field of characteristic zero,
elliptic curve, Somos sequence}

\begin{abstract} We detail the continued fraction expansion of the square
root of the general monic quartic polynomial, noting that each line
of the expansion corresponds to addition of the divisor at infinity.
We analyse the data yielded by the general expansion. In that way we
obtain `elliptic sequences' satisfying Somos relations. I mention
several new results on such sequences. The paper includes a detailed
`reminder exposition' on continued fractions of quadratic
irrationals in function fields.
\end{abstract}

\maketitle
\pagestyle{myheadings}\markboth{{\headlinefont Alf van der
Poorten}}{{\headlinefont Elliptic curves and continued fractions}}

\section{Introduction}

\noindent A delightful `essay' \cite{Za} by Don Zagier explains why the sequence
$(B_h)_{h\in\Z}$, defined by  $B_{-2}=1$, $B_{-1}=1$, $B_0=1$, $B_1=1$, $B_2=1$ and
the recursion 
\begin{equation} \label{eq:Zagier}
B_{h-2}B_{h+3}=B_{h+2} B_{h-1} + B_{h+1}B_{h}\,,
\end{equation}
consists only of integers. Zagier comments that `the proof comes from the
theory of elliptic curves, and can be expressed either in terms of the
denominators of the co-ordinates of the multiples of a particular point on a
particular elliptic curve, or in terms of special values of certain Jacobi theta
functions.'

In the present note I study the \cfe\ of the square root of a quartic
\poly, \emph{inter alia} obtaining sequences generated by recursions such as
\eqref{eq:Zagier}. Here, however, it is clear that I have also constructed the
co-ordinates of the shifted multiples of a point on an elliptic curve and it is it
fairly plain how to relate the surprising integer sequences and the elliptic curves
from which they arise.

A brief reminder exposition on \cf s in quadratic function fields appears as
\S\ref{s:Rappels}, starting at page~\pageref{s:Rappels} below.

It turns out that many of my results related to sequences \`a la \eqref{eq:Zagier}
are contained in the recent thesis \cite{Sw} of Christine Swart; for extended
comment see \S\ref{s:Somos}. Where appropriate, I have added  remarks throughout
the paper. 
Michael Somos, see~\cite{Somos}, had \emph{inter~alia} asked for the
inner meaning of the behaviour of the sequences $(B_h)$, above, and of $(C_h)$
defined by
$C_{h-2}C_{h+2}=C_{h-1}C_{h+1}+C_{h}^2$ and $C_{-2}=1$, $C_{-1}=1$, $C_{0}=1$,
$C_{1}=1$: in the terminology of~\cite{Sw}, of the sequences Somos$(5)$ and
Somos$(4)$. More generally, of course, one may both vary the initial values and
coefficients and generalise the `gap' to $2m$ or $2m+1$ by studying Somos~$2m$,
respectively Somos~$2m+1$, namely sequences satisfying the respective recursions
$$
D_{h-m}D_{h+m}=\sum_{i=1}^m \kappa_iD_{h-m+i}D_{h+m-i} \text{\ or\ }
D_{h-m}D_{h+m+1}=\sum_{i=1}^m \kappa_iD_{h-m+i}D_{h+m-i}\,.
$$
I show in passing, a footnote on page~\pageref{page:satisfying}, that a Somos~$4$
is always a Somos~$6$, while Theorem~\ref{pr:2} points out it is always a
Somos~$5$. After seeing~\cite{Sw}, I added a somewhat painful proof,
Theorem~\ref{th:Somos8} on page~\pageref{ss:Somos8}, that it also always is
a Somos~$8$. For example, Somos$(4)$ satisfies all of
\begin{align*}
C_{h-3}C_{h+3}&=C_{h-1}C_{h+1}+5C_{h}^2\,,\\
C_{h-2}C_{h+3}&=-C_{h-1}C_{h+2}+5C_{h}C_{h+1}\,,\\
C_{h-4}C_{h+4}&=25C_{h-1}C_{h+1}-4C_{h}^2\,.
\end{align*}
In the light of such results one can be confident, see \S\ref{ss:singular},
page~\pageref{ss:singular}, that in general if
$(A_h)$ is a Somos~$4$, equivalently an elliptic sequence satisfying 
$A_{h-2}A_{h+2}=W_2^2A_{h-1}A_{h+1}-W_1W_3A_{h}^2$,
then for all $m$ both
$$
W_1W_2A_{h-m}A_{h+m+1}=W_mW_{m+1}A_{h-1}A_{h+2}-W_{m-1}W_{m+2}A_{h}A_{h+1}
$$
and 
$$
W_1^2A_{h-m}A_{h+m}=W_m^2A_{h-1}A_{h+1}-W_{m-1}W_{m+1}A_{h}^2\,.
$$

\section{Continued Fraction Expansion of the Square Root of a Quartic}
\label{s:cfe}

\noindent We suppose the base field $\F$ is not of characteristic~$2$ because that
case requires nontrivial changes throughout the exposition and not of
characteristic~$3$ because that requires some trivial changes to parts of the
exposition. We study the \cfe\ of a quartic \poly\ $D\in\F[X]$; where $D$ is not a
square. Set
\begin{equation} \label{eq:C}
\mathcal C: Y^2=D(X):=(X^2+f)^2+4v(X-w),
\end{equation}
and for brevity write $A=X^2+f$ and $R=v(X-w)$. For 
$h=0$,
$1$,
$2$,
$\ldots\,$  we denote the \cq s of $Y_0$ by
\begin{equation} \label{eq:notation}Y_h=(Y+A+2e_h)/v_h(X-w_h)\,,
\end{equation}
noting that the $Y_h$ all are reduced, namely $\deg Y_h\gt0$ but $\deg
\conj Y_h\lt0$. The upshot is that the
\cfe\ of
$Y_0$ has typical line, line~$h\,$:
$$
\frac{Y+A+2e_h}{v_h(X-w_h)}=\frac{2(X+w_h)}{v_h} - 
\frac{\conj Y+A+2e_{h+1}}{v_h(X-w_h)}\,.
$$
Thus evident recursion formulas, see \eqref{eq:quadratic} at
page~\pageref{eq:quadratic}, yield
\begin{equation} \label{eq:e} 
f+e_h+e_{h+1}=-w_h^2
\end{equation}
and $-v_hv_{h+1}(X-w_h)(X-w_{h+1})=(Y+A+2e_{h+1})(\conj Y+A+2e_{h+1})$. Hence
\begin{equation} \label{eq:Q} 
v_hv_{h+1}(X-w_h)(X-w_{h+1})=-4(X^2+f+e_{h+1})e_{h+1}+4v(X-w).
\end{equation}
Equating coefficients in \eqref{eq:Q}, and then dividing by $-4e_{h+1}$, we get
\begin{align} \label{eq:X2} 
-4e_{h+1}&=v_hv_{h+1}\,;
\tag*{$X^2$:} \\
 \label{eq:X1}
v/e_{h+1}&=w_h+w_{h+1}\,;
\tag*{$X^1$:} 
\\ \label{eq:X0}
f+e_{h+1}+vw/e_{h+1}&=w_hw_{h+1}\,.
\tag*{$X^0$:} \end{align}
The five displayed equations immediately above readily lead by several routes to
\begin{equation}\label{eq:critical} 
e_he_{h+1}=v(w-w_h)\,.
\end{equation}
For example, apply the remainder theorem to the right hand side of \eqref{eq:Q}
after noting it is divisible by $X-w_h$, and recall \eqref{eq:e}.
\begin{proposition}[Adams and Razar \cite{AR}]\label{pr:1} Denote the two points
at infinity on the elliptic curve ~\eqref{eq:C} by
$S$ and
$O$, with $O$ the zero of its group law. The points
$M_{h+1}:=(w_h,e_h-e_{h+1})$ all lie on~$\mathcal C$. Set
$M_1=M$, and
$M_{h+1}=:M+S_h$. Then $S_h=hS$.
\end{proposition}
\begin{proof} The points $M_{h+1}$ lie on the curve $\mathcal C:Y^2=D(X)$ because
\begin{multline*}
(e_h-e_{h+1})^2-(w_h^2+f)^2
=\bigl(e_h-(e_{h+1}+w_h^2+f)\bigr)\bigl((e_{h}+w_h^2+f)-e_{h+1}\bigr)\\
=-4e_{h}e_{h+1}=4v(w_h-w)\,.
\end{multline*}
The birational transformations
\begin{equation}\label{eq:XY}
X=\bigl(V-v\bigr)\big/U,\qquad Y=2U-(X^2+f)\,;
\end{equation}
conversely,
\begin{equation}\label{eq:UV}
2U= Y+X^2+f\,,\qquad
2V= XY+X^3+fX+2v\,,
\end{equation}
move the point $S$ to $(0,0)$, leave $O$ at infinity, and change the
quartic model to a Weierstrass\ model
\begin{equation} \label{eq:Weierstrass}
\mathcal W: V^2-vV=U^3-fU^2+vwU\,.
\end{equation}
Specifically, one sees that $U(M_{h+1})=-e_{h+1}$, and
$V(M_{h+1})=v-w_he_{h+1}$. We also note that
$U(-M_{h+1})=-e_{h+1}$, $V(-M_{h+1})=w_he_{h+1}$.

To check $S+(M+S_{h-1})=M+S_{h}$ on $\mathcal W$ it suffices for us to show
that the three points $(0,0)$, $(-e_{h}, v-w_{h-1}e_{h})$, and  $(-e_{h+1},
w_{h}e_{h+1})$ lie on a straight line. But that is 
$(v-w_{h-1}e_{h})/e_{h}=w_{h}$. So $w_{h-1}+w_{h}=v/e_{h}$ proves the
claim.
\end{proof}
\section{Elliptic sequences}
\begin{theorem} \label{pr:2}
Let $(A_h)$ be the sequence defined by the `initial' values
$A_0$, $A_1$ and the recursive definition
\begin{equation} \label{eq:elkies}
A_{h-1}A_{h+1}=e_hA_h^2\,.
\end{equation}
Then, given $A_0$, $A_1$, $A_2$, $A_3$, $A_4$ satisfying \eqref{eq:elkies}, the
recursive definition
\begin{equation} \label{eq:even}
A_{h-2}A_{h+2}=v^2A_{h-1}A_{h+1}+v^2(f+w^2)A_h^2
\end{equation}
defines the same sequence as does \eqref{eq:elkies}. Just so,  also
\begin{equation} \label{eq:odd}
A_{h-2}A_{h+3}=-v^2(f+w^2)A_{h-1}A_{h+2}+v^3\bigl(v+2w(f+w^2)\bigr)A_{h}A_{h+1}
\end{equation}
defines that sequence.   
\end{theorem}

\begin{proof} By \eqref{eq:critical} we obtain
\begin{multline*}
e_{h-1}e_h^2e_{h+1}=v^2(w-w_{h-1})(w-w_{h})
\\=v^2(w_{h-1}w_h-w(w_{h-1}+w_h)+w^2)=v^2\left((f+e_h+vw/e_h)-w\cdot
(v/e_h)+w^2\right).
\end{multline*}
Thus
\begin{equation} \label{eq:vital}
e_{h-1}e_h^2e_{h+1}=v^2\left(e_h+(f+w^2)\right).
\end{equation}
 However,
$A_{h-1}A_{h+1}=e_hA_h^2$ entails
$$
A_{h-2}A_{h}A_{h-1}A_{h+1}A_{h}A_{h+2}
=e_{h-1}e_he_{h+1}A_{h-1}^2A_{h}^2A_{h+1}^2\,,
$$
and so $A_{h-2}A_{h+2}=e_{h-1}e_he_{h+1}A_{h-1}A_{h+1}$, which is
\begin{equation} \label{eq:Aeven}
A_{h-2}A_{h+2}=e_{h-1}e_h^2e_{h+1}A_h^2\,.
\end{equation}
On multiplying \eqref{eq:vital} by $A_h^2$ we obtain \eqref{eq:even}.

Similarly \eqref{eq:elkies} yields
$A_{h-1}A_{h+1}A_{h}A_{h+2}=e_he_{h+1}A_h^2A_{h+1}^2$, and so
\begin{equation} \label{eq:odd1}
A_{h-1}A_{h+2}=e_he_{h+1}A_hA_{h+1}\,.
\end{equation}
It follows readily that
\begin{equation} \label{eq:odd2}
A_{h-2}A_{h+3}=e_{h-1}e_h^2e_{h+1}^2e_{h+2}A_hA_{h+1}\,.
\end{equation}
Moreover, \eqref{eq:vital} implies that
\begin{equation*} 
e_{h-1}e_h^3e_{h+1}^3e_{h+2}=v^4\left(e_he_{h+1}+
(f+w^2)(e_h+e_{h+1})+(f+w^2)^2\right)\,.
\end{equation*}
However, by \eqref{eq:e} we know that
$v^2\left(e_h+e_{h+1}+f+w^2\right)=v^2(w^2-w_h^2)$. Here
$v(w-w_h)=e_he_{h+1}$ and $v(w+w_h)=-v(w-w_h)+2vw=-e_he_{h+1}+2vw$. So
\begin{equation} \label{eq:odde}
e_{h-1}e_h^2e_{h+1}^2e_{h+2}=v^2\left(-(f+w^2)e_he_{h+1}+v^2+2vw(f+w^2)\right),
\end{equation}
which immediately allows us to see that also \eqref{eq:odd} yields the sequence
$(A_h)$.
\end{proof}
The extraordinary feature of the identities \eqref{eq:vital} and
\eqref{eq:odde} is their independence of the translation
$M$: thus of the initial data
$v_0$, $w_0$, and
$e_0$. 
\subsection{Two-sided infinite sequences} It is plain that the various
definitions of the sequence $(A_h)$ encourage one to think of it as
bidirectional infinite. Indeed, albeit that one does feel a need to start a
\cfe\ --- so one conventionally begins it at $Y_0$, one is not stopped
from thinking of the tableau listing the lines of the \ex\ as being two-sided
infinite; note the remark at the end of \S~\ref{ss:cfe},
page~\pageref{page:two-sided}. In summary: we may and should view the various
sequences $(e_h)$, $\ldots\,$, defined above, as two-sided infinite sequences.
\label{page:extraordinary}

\subsection{Vanishing} If say $v_k=0$, then line~$k$ of the \cfe\ of
$Y_0$ makes no sense both because the denominator $Q_k(X):=v_k(X-w_k)$ of the
\cq\
$Y_k$ seems to vanish identically and because the alleged \pq\
$a_k:=2(X+w_k)/v_k$ blows up. 

The second difficulty is real. The vanishing of $v_k$ entails a \pq\ blowing up
to higher degree. We deal with vanishing by refusing to look at it. We move the
point of impact of the issue by dismissing most of the data we have obtained,
including the \cf\ tableau, and keep only a part of the sequence $(e_h)$. That
makes the first difficulty moot.\footnote{In any case, the first apparent
difficulty is just an artifact of our notation. If, from the start, we had
written
$Q_h=v_hX+y_h$, as we might well have done at the cost of nasty fractions in our
formulas, we would not have entertained the thought that
$v_k=0$ entails $y_k=0$. Plainly, we must allow $v_k=0$ yet $v_kw_k\ne0$.}

\begin{remark*} There is no loss of generality in
taking
$k=0$. Then, up to an irrelevant normalisation, $Y_0=Y+A$. 
If more than one of the $v_h$ vanish then it is a simple exercise
to confirm that the
\cfe\ of
$Y_0$ necessarily is purely periodic, see the discussion at
page~\pageref{page:quasiperiodic}. If $Y_0$ does not have a periodic \cfe\ then
there is some $h_0$, namely $h_0=0$, so that, for all $h\gt h_0$, line~$h$ of the
\ex\ of
$Y_0$ does make sense. 
\end{remark*}
Except of course when dealing explicitly with periodicity, we suppose in the
sequel that if $v_k=0$ then $k=0$; we refer to this case as the \emph{singular}
case.

\subsection{The singular case} \label{ss:singular}

We remark that in the singular
case the sequence
$(e_h)_{h\ge1}$ defines antisymmetric double-sided sequences $(W_h)$, that is
with $W_{-h}=-W_h$, by
$W_{h-1}W_{h+1}=e_hW_h^2$ and so that, for all integers $h$, $m$, and $n$,
\begin{equation} \label{eq:redundant'}
W_{h-m}W_{h+m}W_n^2+W_{n-h}W_{n+h}W_m^2
+W_{m-n}W_{m+n}W_h^2=0\,.\tag{\ref{eq:redundant}$'$}
\end{equation} 
Actually, one may find it preferable to forego an insistence on
antisymmetry in favour of rewriting \eqref{eq:redundant'} less elegantly as
\begin{equation} \label{eq:redundant}
W_{h-m}W_{h+m}W_n^2=W_{h-n}W_{h+n}W_m^2 - W_{m-n}W_{m+n}W_h^2\,,
\end{equation}
just for $h\ge m\ge n$. 
In any case, \eqref{eq:redundant} seems more dramatic than it is. An
easy exercise confirms that, if $W_1=1$,~\eqref{eq:redundant} is equivalent
to just 
\begin{equation} \label{eq:general}
W_{h-m}W_{h+m}=W_m^2 W_{h-1}W_{h+1} -W_{m-1}W_{m+1}W_h^2
\end{equation}
for all integers $h\ge m$. Indeed, \eqref{eq:general} is just a special
case of \eqref{eq:redundant}. However, given~\eqref{eq:general}, obvious
substitutions in
\eqref{eq:redundant} quickly show one may return from
\eqref{eq:general} to the apparently more general \eqref{eq:redundant}.

But there is a drama here. As already remarked in a near identical situation, the
recurrence relation $W_{h-2}W_{h+2}=W_2^2 W_{h-1}W_{h+1} -W_{1}W_{3}W_h^2$, and
five or so initial values, already suffices to produce $(W_h)$. Thus
\eqref{eq:general} for all $m$ is apparently entailed by its special case $m=2$.

I can show this directly\label{page:Somos6}\footnote{Plainly
$e_{h-2}e_{h-1}^2e_{h}^3e_{h+1}^2e_{h+2}\cdot e_{h}=v^4\bigl(e_{h-1}+(f+w^2)\bigr)
\bigl(e_{h+1}+(f+w^2)\bigr)e_h^2$. Now notice that
$(e_{h-1}e_{h}+e_{h}e_{h+1})e_{h}=v(w-w_{h-1}+w-w_h)e_h=2vwe_h-v^2$ and recall
that
$e_{h-1}e_{h}^2e_{h+1}=v^2\bigl(e_{h}+(f+w^2)\bigr)$. The upshot is a miraculous
cancellation yielding 
$$
e_{h-2}e_{h-1}^2e_{h}^3e_{h+1}^2e_{h+2}\cdot
e_{h}=v^4\bigl((f+w^2)^2e_h^2+v(v+2w(f+w^2))e_h\bigr)
$$
and allowing us to divide by the auxiliary $e_h$. Thus the bottom line is 
$$
A_{h-3}A_{h+3}=v^4\bigl((f+w^2)^2A_{h-1}A_{h+1}+v(v+2w(f+w^2))A_{h}^2\bigr),
$$
which is
$A_{h-3}A_{h+3}=W_3^2A_{h-1}A_{h+1}-W_2W_4A_{h}^2$.}, 
by way
of new relations on the
$e_h$, for
$m=3$. But the case $m=4$ already did not seem worth the
effort. Whatever, my approach
gave me no hint as to how to concoct an inductive argument leading to general
$m$. Plan B, to look it up, fared little better. In her thesis \cite{Shi}, Rachel
Shipsey shyly refers the reader back to Morgan Ward's opus \cite{Wa}; but Ward
does not comment on the matter at all, having \emph{defined} his sequences by
\eqref{eq:general}.  Well, perhaps Ward does comment. The issue is whether
\eqref{eq:general} is coherent: do different $m$ yield the one sequence? Ward
notes that if $\sigma$ is the Weierstra\ss\
$\sigma$-function then a sequence $\bigl(\sigma(hu)/\sigma(u)^{h^2}\,\bigr)$
satisfies
\eqref{eq:general} for all $m$. Whatever, a much more direct argument would be
much more satisfying.\footnote{For additional remarks, and a dissatisfying proof
for the case $m=4$, see
\S\ref{ss:Somos8} at page~\pageref{ss:Somos8}.}\label{page:satisfying}

Proposition~\ref{pr:2} shows that certainly $W_{h-2}W_{h+2}=W_2^2 W_{h-1}W_{h+1}
-W_{1}W_{3}W_h^2$ for $h=1$, $2$, $\ldots\,$, in which case \eqref{eq:general}
apparently follows by arguments in \cite{Wa} and anti-symmetry; 
\eqref{eq:redundant} is then just an easy exercise.

The singular case is initiated by $v_1=4v$, $w_1=w$, $e_1=0$, $e_2=-(f+w^2)$. For
temporary convenience set $x=v/(f+w^2)$. From the original \cfe\ of $Y+A$
or, better, the recursion formulas of page~\pageref{eq:e}, we fairly readily
obtain $v_2=1/x$, $w_2=w-x$, $e_3=-x(x+2w)$, 
$e_4=v(x^2(x+2w)-v)/x^2(x+2w)^2$.

We are now free to choose, say $W_1=1$, $W_2=v$, leading to $W_3=-v^2(f+w^2)$,
$W_4=-v^4\bigl(v+2w(f+w^2)\bigr)$,
$W_5=-v^6\bigl(v(v+2w(f+w^2))-(f+w^2)^3\bigr)$, $\ldots\,$.

That allows us to notice that \eqref{eq:odd} apparently is
$$
W_1W_2A_{h-2}A_{h+3}=W_2W_3A_{h-1}A_{h+2}-W_1W_4A_{h}A_{h+1}
$$
and that \eqref{eq:even} of course is
$$
A_{h-2}A_{h+2}=W_2^2A_{h-1}A_{h+1}-W_1W_3A_{h}^2\,.
$$
Given that also $A_{h-3}A_{h+3}=W_3^2A_{h-1}A_{h+1}-W_2W_4A_{h}^2$, it is of course
tempting to guess that more is true. Certainly, more \emph{is} true in the
special case $(A_h)=(W_h)$, that's the point of the discussion above. Moreover,
the same `more' is true, see for example \cite[Theorem~8.1.2, p.\,191]{Sw}, for
sequence translates: thus $(A_h)=(W_{h+k})$.\label{page:more}
\subsection{Elliptic divisibility sequences}\label{ss:eds}  Recall that in the
singular case and for
$h=1$,
$2$,
$\ldots\,$ the $-e_h$ are in fact the $U$ co-ordinates of the multiples $hS$ of
the point $S=(0,0)$ on the curve 
$V^2-vV=U^3-fU^2+vwU$. 

Suppose we are working in the ring $Z=\Z[f, v, vw]$ of `integers'. If
 $\gcd(v,vw)=1$, so the exact denominator of the `rational' $w$ is $v$,
then our choices $W_1=1$, $W_2=v$ lead the definition
$W_{h-1}W_{h+1}=e_hW_h^2$ to be such that $W_h^2$ is always the exact
denominator of the `rational' $e_h$. It is this that is shown in detail 
by Rachel Shipsey \cite{Shi}. In particular it follows that $(W_h)$ is an
elliptic divisibility sequence as described by Ward \cite{Wa}. A convenient
recent introductory reference is Chapter~10 of the book \cite{EPSW}.

Set $hS=(U_h/W_h^2, V_h/W_h^3)$, thus
defining sequences $(U_h)$, $(V_h)$, and $(W_h)$ of integers, with $W_h$ chosen
minimally. Shipsey notices, provided that indeed $\gcd(v,vw)=1$, that
$\gcd(U_h,V_h)=W_{h-1}$ and $W_{h-1}W_{h+1}=-U_h$. Here, I have used this last
fact to define the sequence $(W_h)$. 

Starting, in effect, from the definition \eqref{eq:redundant}, Ward
\cite{Wa} shows that with $W_0=0$, $W_1=1$, and $W_2\Div W_4$, the sequence
$(W_h)$ is a divisibility sequence; that is, if $a\Div b$ then $W_a\Div W_b$. A
little more is true. If also $\gcd(W_3,W_4)=1$ then in fact
$\gcd(W_a,W_b)=W_{\gcd(a,b)}$. On the other hand, a prime dividing both $W_3$ and
$W_4$ divides $W_h$ for all $h\ge3$.

\subsection{Periodicity}  Suppose now that the sequence $(W_h)$ is periodic.  From
the \cfe\ and, say, \cite{163}, we find that $v=0$ (but $w'=vw\ne0$ if our curve
is to be elliptic) is the case of the \cf\ having  quasi-period $r=1$ and
the divisor at infinity on the curve having torsion
$m=2$. Just so, $f+w^2=0$, thus $W_3=0$, signals $r=2$ and $m=3$, and $x+2w=0$,
or
$W_4=0$, is $r=3$ and $m=4$. And so on; for more see \cite{163}. In summary, 
$m\gt0$ is minimal with $W_m=0$ if and only if the \cfe\ of, say,
$Z=\frac12(Y+A)$ has a minimal quasi-period of length $r=m-1$. Note that by
observations summarised at \S\ref{ss:units}, page~\pageref{ss:units} below, the
minimal period of the \cfe\ of $Z$ then is $m-1$ or $2m-2$ if $m$ is even, and
$m-1$ if $m$ is odd. This holds\footnote{This fact is not well understood. For
example, over the finite field $\F_p$ of $p$ elements, $p$ an odd prime, it is
plain and well known that a quasi-period of length $r$ entails there is a period of
length
$rt$ for some divisor $t$ of $p-1$. Indeed, if we are expanding a general
quadratic irrational function, then this is the best one can assert. If, however,
the function has integral --- thus \poly\ --- trace, \emph{a fortiori} if it is a
quadratic irrational  integer element, then if its quasi-period is of even
length it must already be a period, and if it of odd length $r$, then at any rate
$2r$ is a period~length.} for any base field. Over
$\Q$, recall~\cite{Ma}, the only possibilities for $m$ are
$2$,
$3$, $\ldots\,$, $10$, or $12$; by happenstance, the cases period $9$
or
$11$ do not occur (see \cite{149} or \cite{163}).

However, Ward~\cite{Wa} shows and one fairly
readily confirms that precisely the periods $1$, $2$, $3$, $4$, $5$, $6$, $8$, or
$10$ are possible for an integral elliptic divisibility sequence defined by
\eqref{eq:general}. A convenient summary and explanation is given by Chistine
Swart~\cite[\S4.6]{Sw}.

\begin{aside}\label{a:aside} It has been suggested in my hearing that `Mathematics
is the study of degeneracy', so the following warrants careful consideration. In
the singular case we have
$e_1=0$ and then the recursion 
$e_{h-1}e_h^2e_{h+1}=v^2\bigl(e_h+(f+w^2)\bigr)$ and $e_2=-(f+w^2)$ forces 
$e_0\cdot 0^2=-v^2$ in the case $h=1$. In a context in which $v=0$ and $w'=vw\ne0$
passes without comment this is no great matter. However, we must also cope with 
$W_{-1}W_{1}=e_0W_0^2$, and with
both $W_{-2}W_{0}=e_{-1}W_{-1}^2$ and $W_{-3}W_{-1}=e_{-2}^{}W_{-2}^2$ and it now
is more difficult to reconcile the two-sided \cfe\ and the two-sided
anti-symmetric singular elliptic sequence. The periodic case seems particularly
problematic. Noticeably, Christine Swart~\cite{Sw} declares her elliptic sequences
as undefined beyond a $0$; I've chosen to be vague. 

\end{aside}


\section{Examples} 

\subsection{}\label{ss:1}  Consider the curve $\mathcal
C:Y^2=(X^2-29)^2-4\cdot48(X+5)$; here a corresponding cubic model is $\mathcal E:
V^2+48V=U^3+29U^2+240U$. Set
$A=X^2-29$. The first several preceding and succeeding steps in the \cfe\ of
$Y_0=(Y+A+16)/8(X+3)$ are\footnote{Here my choice of $v_0=8$ is arbitrary but not
at random.}
\begin{align}
\frac{Y+A+18}{16(X+2)/3}&=\frac{X-2}{8/3}-\frac{\conj Y+A+32}{16(X+2)/3}
\tag*{line $\conj3:$}\\
\frac{Y+A+32}{12(X+1)}&=\frac{X-1}{6}-\frac{\conj Y+A+24}{12(X+1)}
\tag*{line $\conj2:$}\\
\frac{Y+A+24}{4(X+3)}&=\frac{X-3}{2}-\frac{\conj Y+A+16}{4(X+3)}
\tag*{line $\conj1:$}\\
\frac{Y+A+16}{8(X+3)}&=\frac{X-3}{4}-\frac{\conj Y+A+24}{8(X+3)}
\tag*{line $0:$}\\
\frac{Y+A+24}{6(X+1)}&=\frac{X-1}{3}-\frac{\conj Y+A+32}{6(X+1)}
\tag*{line $1:$}\\
\frac{Y+A+32}{32(X+2)/3}&=\frac{X-2}{16/3}-\frac{\conj Y+A+18}{32(X+2)/3}
\tag*{line $2:$}\\
\frac{Y+A+18}{9(3X+10)/8}&=\;\cdots
\tag*{} 
\end{align}
where elegance has suggested we write `line $\conj h$' as short for `line $-h$'.

The feature motivating this example is the six integral points $(-2,\pm7)$,
$(-1,\pm4)$, and $(-3,\pm4)$ on $\mathcal C$. With $M_{\mathcal C}=(-3,4)$ and
$S_{\mathcal C}$ the `other' point at infinity these are in fact the six points
$M_{\mathcal C}+hS_{\mathcal C}$ for $h=-3$,
$-2$, $-1$,
$0$,
$1$, and
$2$.

Correspondingly, on $\mathcal E$ we have the integral points $M+2S=(-16,-32)$ and
$M-2S=(-16,-16)$, $M-S=(-12,-36)$ and $M+S=(-12,-12)$;  here
$M=M_{\mathcal E}=(-8,-24)$; $S=S_{\mathcal
E}=(0,0)$. Of course $\mathcal E$ is not minimal; nor, for that matter was
$\mathcal C$.  In fact the replacements $X\leftarrow 2X+1$, $Y\leftarrow 4Y$  yield
\begin{equation}\label{eq:min}
Y^2=(X^2+X-7)^2-4\cdot6(X+3)\,,
\end{equation}
correctly suggesting we need a more general treatment than that presented in the
discussion above. It turns out to be enough  for present purposes to replace
$e_h\leftarrow 4e_h$ obtaining
$$\dots\,,\; e_{-3}=\tfrac{9}{4},\; e_{-2}=4,\; e_{-1}=3,\; e_{0}=2,\;
e_{1}=3,\; e_{2}=4,\; e_{3}=\tfrac{9}{4},\;\ldots\;.
$$
Then $A_0=1$, $A_1=1$ and 
$$A_{h-1}A_{h+1}=e_hA_h^2
$$
yields the sequence $\ldots\,$, $A_{-4}=2^53^5$, $A_{-3}=2^53^2$, $A_{-2}=2^33$,
$A_{-1}=2$, $A_{0}=1$, $A_{1}=1$, $A_{2}=3$, $A_{3}=2^23^2$, $A_{4}=2^23^5$,
$\ldots\,$.
Notice that we're hit for six\footnote{My remark is
guided by knowing that $V^2+UV+6V=U^3+7U^2+12U$ is a minimal model for $\mathcal
E$, and noticing that $\gcd(6, 12)=6$.} by increasingly high powers of primes
dividing
$6$ appearing as factors of the $A_h$. However, we know that \eqref{eq:odd}
derives from \eqref{eq:odde}. With the original $e_h$s divided by $4$ that yields
\begin{equation*} 
A_{h-2}A_{h+3}=6^2A_{h-1}A_{h+2}+6^3A_{h}A_{h+1}\,.
\end{equation*}
Remarkably, one may remove the effect of the $6$ by renormalising to a sequence
$(B_h)$ of integers satisfying 
$$
B_{h-2}B_{h+3}=B_{h-1}B_{h+2}+B_{h}B_{h+1}\,.
$$
Specifically, $\ldots\,$, $B_{-4}=3$, $B_{-3}=2$, $B_{-2}=1$,
$B_{-1}=1$, $B_{0}=1$, $B_{1}=1$, $B_{2}=1$, $B_{3}=2$, $B_{4}=3$, $B_{5}=5$,
$B_{6}=11$, $B_{7}=37$, $B_{8}=83$, $\ldots\,$, and the sequence is symmetric
about $B={0}$. Interestingly, the choice of each $B_h$ as a divisor of $A_h$ is
forced, in the present case by the data $A_0=A_1=1$ and the decision that the
coefficient of $B_hB_{h+1}$ be $1$.  Of course it is straightforward to verify
that $A_{h-2}A_{h+3}$ is always divisible by $6^3$ and $A_{h-1}A_{h+2}$ always by
$6$. For a different treatment see \S\ref{ss:4} below.

\subsection{} 

Take $v=\pm1$ and $f+w^2=1$. Thus $e_{h-1}e_h^2e_{h+1}=e_h+1$
and so $e_0=1$, $e_1=1$ yields the sequence $\ldots\,$, $2$, $1$, $1$, $2$,
$3/4$, ${14}/{9}$, $\ldots\,$, of values of $e_h$. As explained above, with $C_0=1$
and
$C_1=1$, the definition $C_{h-1}C_{h+1}=e_hC_h^2$ yields the symmetric
sequence
$\ldots\,$,
$2$, $1$,
$1$, $1$, $1$,
$2$, $3$, $7$, $23$, $59$, $\ldots\,$, of values of $C_h$ satisfying the recursion
$$C_{h-2}C_{h+2}=C_{h-1}C_{h+1}+C_h^2\,.
$$
Plainly, one can get four consecutive values $1$ in a sequence $(C_h)$ as just
defined only by having two consecutive values $1$ in the corresponding sequence
$(e_h)$.

Set $Y^2=A^2+4v(X-w)$, where $A=X^2+f$. With $Z=\frac12(Y+A)$, we have $Z\ov
Z=-v(X-w)$ and $Z+\ov Z=A$. Thus the consecutive lines
\begin{align*}
\frac{Z+1}{X-b}&=\phantom{-}(X+b) - \frac{\ov Z+1}{X-b}\\
\frac{Z+1}{-(X-c)}&=-(X+c) - \frac{\ov Z+2}{-(X-c)}
\end{align*}
entail 
$$f+2=-b^2\,,\quad f+3=-c^2\,,\quad \text{and}\quad b+c=v\,,\quad bc=f+vw+1\,,
$$
which is $v=\pm1$, $b=\pm1$, $c=0$, $f=-3$, and $w=\pm2$. Up to $X\gets
-X$, the sequence $(C_h)$ is given by the curve $\mathcal C: Y^2=(X^2-3)^2+4(X-2)$
and its points $M_{\mathcal C}+hS_{\mathcal C}$, $M_{\mathcal C}=(1,0)$,
$S_{\mathcal C}$ the `other point' at infinity; equivalently by 
$$\mathcal E:V^2-V=U^3+3U^2+2U \quad\text{with $M=(-1,1)$,
$S=(0,0)$}.
$$
Indeed, $M+S=(-1,0)$, $M+2S=(-2,1)$, $M+2S=(-3/4,3/8)$, $\ldots\,$.  Note that it
is impossible to have three consecutive values
$1$ in the sequence
$(e_h)$ if also $v=\pm1$, except for trivial periodic cases, so the
hoo-ha of the example at
\S\ref{ss:1} above is in a sense unavoidable.

\subsection{Remarks}\label{ss:3} The two examples get  a rather woolly
treatment at \cite{Sp} and its preceding discussion, see
\cite{Somos} for context. Before seeing \cite{Sw} I had also remarked that ``the
observation that a twist
$V^2-vV=dU^3-fU^2+vwU$ becomes $V^2-dvV=U^3-fU^2+dvw$ by $U\gets dU$, $V\gets dV$
allows one to presume
$v=\pm1$. A suitable choice of $e_0$, $e_1$ and $A_0$, $A_1$ should now allow one
to duplicate the result 
 claimed in \cite{Sp} in
 somewhat less brutal form." Namely, one wishes to obtain elliptic
curves yielding a sequence $(A_h)$ with nominated $A_{-1}$, $A_0$, $A_1$, $A_2$
and such that $A_{h-2}A_{h+2}=\kappa A_{h-2}A_{h-2}+\lambda A_{h-2}A_{h-2}$;
only the cases $\kappa$ not a square are at issue. In fact, this issue is dealt
with by Christine Swart at
\cite[p.\,153\emph{ff}]{Sw} in more straightforward fashion than I had in mind. In
very brief, if
$A_{h-1}A_{h+1}=e_hA_h^2$, then
$$B_h=\kappa^{\frac12h(h+1)}A_h \quad\text{entails}\quad B_{h-1}B_{h+1}=
\kappa e_hB_h^2\,,
$$
and so $B_{h-2}B_{h+2}=(1/\kappa)^2 B_{h-2}B_{h-2}+(\lambda/\kappa^4)
B_{h-2}B_{h-2}$.
\subsection{Reprise}\label{ss:4} It seems appropriate to return to the example of
\S\ref{ss:1} so as to \emph{discover} the elliptic curve giving rise to
$(B_h)=(\ldots, 1, 1, 1, 1, 1, \ldots)$, given that
$B_{h-2}B_{h+3}=B_{h-1}B_{h+2}+B_{h}B_{h+1}$. Recall we expect the squares of the
integers $B_h$ to be the precise denominators of the points $M+hS$ on the minimal
Weierstra\ss\ model ${\mathcal W}$ of the curve; here
$M$ is some point on that model and $S=(0,0)$.

Suppose $e_{-2}$, $e_{-1}$, $e_{0}$, $e_{1}$, $e_{2}$ supply the five
integer co-ordinates yielding $B_{-2}$, $B_{-1}$, $B_{0}$, $B_{1}$, $B_{2}$. Of
course no more than two of these $e_i$ can be $1$ so we must have
$$e_0B_{-1}B_1=e_0B_0^2\,,\quad e_1B_{0}B_2=e_1B_1^2\,,\quad
\tfrac12 e_2B_{1}B_3=e_2B_2^2\,,
$$
since of course the recursion for $(B_h)$ entails $B_3=2$. Suppose in general that
\begin{equation*} 
c_hB_{h-1}B_{h+1}=e_hB_h^2\,. 
\end{equation*}
Then the identity \eqref{eq:odde}, 
\begin{equation*}
e_{h-1}e_h^2e_{h+1}^2e_{h+2}=v^2\left(-(f+w^2)e_he_{h+1}+v^2+2vw(f+w^2)\right),
\end{equation*}
and $B_{h-2}B_{h+3}=B_{h-1}B_{h+2}+B_{h}B_{h+1}$ entail  
$$
c_{h-1}c_h^2c_{h+1}^2c_{h+2}=-v^2(f+w^2)c_hc_{h+1}=v^3\bigl(v+2w(f+w^2)\bigr).
$$
Thus $c_hc_{h+1}=kv$, say, is independent of $h$ and we have
\begin{equation*} 
k^2=-(f+w^2)\quad\text{and}\quad k^3+2wk^2-v=0\,.
\end{equation*}
Note that if $f+w^2$, or $2w$ and $v$, are integers, also $k$ must be an integer.
Also,
\begin{equation}\label{eq:vk}
e_0e_1=vk\quad\text{and}\quad e_1e_2=2vk\,.
\end{equation}

\begin{Remark*} However, 
$e_{h-1}e_h^2e_{h+1}=v^2e_h+v^2(f+w^2)$ implies 
\begin{equation*}
k^2B_{h-2}B_{h+2}=c_hB_{h-1}B_{h+1}+(f+w^2)B_h^2\,,
\end{equation*}
without the coefficients necessarily being independent of $h$.
In particular, $k^2=-(f+w^2)$ entails  $c_0=e_0=2k^2$ and $c_1=e_1=3k^2$. 
\end{Remark*}

On the
other hand, the identity \eqref{eq:critical} now reports that $k=w-w_0$ and
$2k=w-w_1$. By
\eqref{eq:Q} we then have
\begin{equation*}
w_0+w_1=v/e_1=2w-3k
 \end{equation*}
whilst by \eqref{eq:e} we see that $f+e_0+e_1=-(w-k)^2$, $f+e_1+e_2=-(w-2k)^2$,
so, recalling that $e_2=2e_0$, also
$e_0=(2w-3k)k$. Hence 
\begin{equation}\label{eq:w}
(2w-3k)k=2k^2 \quad\text{and so}\quad 2w=5k\,.  \end{equation}
In summary, we then can quickly conclude that also
\begin{equation}\label{eq:f}
v=6k^3\,,\quad 4f=-29k^2\,, \quad\text{and}\quad 2w_0=3k\,.  \end{equation}
 
The normalisation $k=-2$ retrieves the \cfe\ at \S\ref{ss:1} on
page~\pageref{ss:1}. As shown at \S\ref{s:comments} on
page~\pageref{s:comments} the corresponding minimal Weierstra\ss\ model is
$V^2+UV+6V=U^3+7U^2+12U$; and  $M=(-2,-2)$ is a point of order two. 

\section{Somos Sequences}\label{s:Somos}

\subsection{Christine Swart's thesis \cite{Sw}}\label{ss:Chris} Much of the
work reported by me here is {\it sui generis\/} with my original object having
been little more than to make explicit the ideas of Adams and Razar~\cite{AR} and
to rediscover elliptic torsion surfaces by a seemingly new method, see
\cite{163} and its references. Eventually, I learned of Michael Somos's sequences,
see~\cite{Somos}, and realised how they arise from my data. I had sort of heard of
Christine Swart's work from Nelson Stephens in 2003 but, fortunately as it turned
out\footnote{It was fun to puzzle out not just the answers to questions, but also
to attempt to guess what the questions ought to be.}, did not have access
to her thesis~\cite{Sw} until very recently when this paper was already essentially
complete; see~\cite{169}.

Christine Swart's discussion of the interrelationship between elliptic sequences
and elliptic curves is more detailed and complete than mine. Among many other
things, she is careful to recognise that the formulas neither know nor care
whether the given elliptic curve is in fact elliptic: thus, for example, also my
quartics may have multiple zeros. If so, extra comment --- mostly, quite
straightforward ---  is required at a number of points above: but is neglected by
me. Further, Christine Swart reports \emph{inter alia} that Nelson Stephens
(personal communication to her) had noticed identities equivalent to
\eqref{eq:vital} and
\eqref{eq:even}, both at page~\pageref{eq:even}; these are her Theorems 3.5.1
and 7.1.2
\cite[p.\,29 and p.\,153]{Sw}.

\subsection{Much more satisfying?}\label{ss:Somos8} I can in fact show that an
elliptic sequence $(A_h)$ also is given by 
$A_{h-4}A_{h+4}=W_4^2A_{h-1}A_{h+1}-W_3W_5A_h^2$. However, I consider my argument
as typical of the sort of thing that gives mathematics its bad name: and
regret to have to admit that this sort of nonsense seemingly does generalise to
proving that a Somos~4 also always is a three-term Somos~$(4+n)$.

Specifically, we know that 
\begin{align}
A_{h-2}A_{h+2}&=W_2^2A_{h-1}A_{h+1}-W_1W_3A_h^2\,, \quad\text{by
definition;}\label{eq:Somos4}\\
A_{h-3}A_{h+3}&=W_3^2A_{h-1}A_{h+1}-W_2W_4A_h^2\,,\quad\text{see footnote
page~\pageref{page:satisfying},}\label{eq:Somos6}\\
\intertext{and intend to show that}
A_{h-4}A_{h+4}&=W_4^2A_{h-1}A_{h+1}-W_3W_5A_h^2\label{eq:Somos8}\,.
\end{align}
Notice that, in particular,
$W_1W_5=W_2^2W_2W_4-W_1W_3W_3^2$; so because $W_1=1$,  $W_3^4=W_2^3W_3W_4-W_3W_5$.
Now observe that 
\begin{equation*}
A_{h-4}A_{h+2}A_{h-2}A_{h+4}
=(W_3^2A_{h-2}A_{h}-W_2W_4A_{h-1}^2)(W_3^2A_{h}A_{h+2}-W_2W_4A_{h+1}^2)\,.
\end{equation*}
In the product on the right, the first term is 
$$
(W_2^3W_3W_4A_h^2-W_3W_5A_h^2)A_{h-2}A_{h+2}
$$
and half of it contributes half of \eqref{eq:Somos8}. Similarly, half the final
term of the product, thus of
$$W_4^2A_{h-1}A_{h+1}\cdot W_2^2A_{h-1}A_{h+1}=
W_4^2A_{h-1}A_{h+1}(A_{h-2}A_{h+2}+W_3A_h^2)\,,
$$
provides the other half of \eqref{eq:Somos8}. Thus it's ugly but true that we have
proved that \eqref{eq:Somos8} holds if and only if $W_3W_4=0$ or
\begin{equation*}
W_2^3A_{h-2}A_h^2A_{h+2}-W_2W_3(A_{h-2}A_hA_{h+1}^2
+A_{h-1}^2A_hA_{h+2})+W_4A_{h-1}A_h^2A_{h+1}=0\,.
\end{equation*}
I now compound this brutality by fiercely replacing the two occurrences of
$A_{h-2}$ by the evident relation
$$
A_{h-2}=(W_2^2A_{h-1}A_{h+1}-W_1W_3A_h^2)/A_{h+2}\,.
$$
That necessitates our then multiplying by $A_{h+2}$. Fortunately, we can
compensate for this cruelty by dividing by $A_h$. We are left with needing to show
that
\begin{multline}\label{eq:!}
W_2^5A_{h-1}A_hA_{h+1}A_{h+2}-W_3^2W_3A_h^3A_{h+2}-W_2^3W_3A_{h-1}A_{h+1}^3\\
-W_2W_3A_{h-1}^2A_{h+2}^2+W_2W_3^2A_h^2A_{h+1}^2
+W_4A_{h-1}A_hA_{h+1}A_{h+2}=0\,.
\end{multline}
It's now natural to despair, and to start looking for a Plan B. However, one
might notice, on page~\pageref{page:satisfying}, that
$W_4=-v^4\bigl(v+2w(f+w^2)\bigr)$; and
$W_2=v$. Moreover
$W_3=-v^2(f+w^2)$. Thus, conveniently, 
\begin{equation*}
W_2^5+W_4=-2v^4w(f+w^2)=W_2^2W_3\,.
\end{equation*}
Hence, just as our result is trivial if $W_3W_4=0$, so also it is trivial if
$W_2=0$. All this is a sign that we may not as yet have made an error. We may
divide
\eqref{eq:!} by 
$W_2W_3$. Better yet, let's also divide by $A_h^2A_{h+1}^2$ by using the
definitions
$$
A_{h-1}A_{h+1}=e_hA_h^2\,, \quad\text{whence also}\quad
A_{h-1}A_{h+2}=e_he_{h+1}A_hA_{h+1}\,. 
$$
Then all that remains  is a confirmation that
\begin{equation}\label{eq:quartic}
2vwe_he_{h+1}-v^2(e_h+e_{h+1})-e_h^2e_{h+1}^2-v^2(f+w^2)=0\,.
\end{equation}
However \eqref{eq:vital}, page~\pageref{eq:vital}, is $e_he_{h+1}=v(w-w_h)$, while
$e_h+e_{h+1}=-f-w_h^2$ is \eqref{eq:e}, page~\pageref{eq:e}. Astonishingly, the
claim \eqref{eq:quartic} follows immediately.

\begin{theorem} If $(A_h)$ is a Somos~$4$ then it is a Somos~$8$ of the shape
$$
A_{h-4}A_{h+4}=\kappa A_{h-1}A_{h+1}+\lambda A_h^2\,.
$$
\label{th:Somos8}
\end{theorem}
\begin{proof} Given the argument above, it suffices to note that Christine Swart
\cite[p.\,153\emph{ff}]{Sw} shows that any Somos~4 is equivalent to an elliptic
sequence.
\end{proof}

\section{Rappels}\label{s:Rappels}

\subsection{Continued fraction expansion of a quadratic irrational}
\label{ss:cfe}
Let $Y=Y(X)$ be a quadratic irrational integral element of the field
$\F((X^{-1}))$ of Laurent series
\begin{equation} \label{eq:Laurent}
\sum_{h=-d}^\infty f_{-h}X^{-h}, \quad\text{some $d\in\Z$}\quad
\end{equation}
defined over some given base field $\F$; that is, there are \poly s
$T$ and $D$ defined over $\F$ so that 
\begin{equation} 
\label{eq:eq}Y^2=T(X)Y+D(X)\,.
\end{equation}
Plainly, by translating $Y$ by a \poly\ if necessary, we may
suppose that $\deg D\ge 2\deg T+2$, with $\deg D=2g+2$, say, and
$\deg T\le g$; then $\deg Y=g+1$. Recall here that a Laurent series
\eqref{eq:Laurent} with $f_{d}\ne0$ has degree~$d$.

Set 
$Y_0=(Y+P_0)/Q_0$ where $P_0$ and $Q_0$ are \poly s so
that
$Q_0$ divides the norm
$(Y+P_0)(\conj Y+P_0)$; notice here that an $\F[X]$-module
$\id{Q, Y+P}$ is an ideal in $\F[X,Y]$ if and only if $Q\Div
(Y+P)(\conj Y+P)$. 

Further, suppose that
$\deg Y_0\gt0$ and
$\deg \conj Y_0\lt 0$; that is, $Y_0$ is \emph{reduced}.
Then the
\cfe\ of $Y_0$ is given by a sequence of lines, of which the
$h$-th is
\begin{equation} \label{eq:lineh}
Y_h:=(Y+P_h)/Q_h=a_h-(\overline Y+P_{h+1})/Q_h\,;
\quad\text{in
brief}\quad Y_h=a_h-\conj B_h 
\,.
\end{equation}
Here the \poly\ $a_h$ is a \emph{\pq}, and the next \emph{complete
quotient}
$Y_{h+1}$ is the reciprocal of the preceding \emph{remainder}
$-(\overline Y+P_{h+1})/Q_h$\,.  Plainly the sequences of \poly s
$(P_h)$ and $(Q_h)$ are given by the recursion formulas
\begin{equation}
\label{eq:quadratic}P_h+P_{h+1}+(Y+\ov Y\,)=a_hQ_h
\ \text{and}\ 
Y\conj Y+(Y+\conj Y)P_{h+1}+P_{h+1}^2=-Q_hQ_{h+1}\,.
\end{equation}
It is easy to see by induction on $h$ that
$Q_h$ divides the norm
$(Y+P_h)(\conj Y+P_h)$.

We observe also that we have a conjugate expansion with $h$-th line
\begin{equation} \label{eq:conjlineh}
B_h:=(Y+P_{h+1})/Q_h=a_h-(\overline Y+P_{h})/Q_h\,,
\quad\text{that is,}\quad 
B_h=a_h-\conj Y_h\,.
\end{equation}
Note that the next line of this expansion is the
conjugate of the previous line of its conjugate \ex\,: conjugation
reverses a \cf\ tableau. Because the conjugate of line~$0$ is the last line of
its tableau we can extend the \ex\ forming the conjugate tableau, leading to
lines $h=1$, $2$, $\ldots\,$
$$(Y+P_{-h+1})/Q_{-h}=a_{-h}-(Y+P_{-h})/Q_{-h}\,; \quad\text{that is,}\quad
B_{-h}=a_{-h}-\conj Y_{-h}\,.
$$
Plainly the original \cf\ tableau also is two-sided infinite and our thinking of
it as `starting' at $Y_0$ is just convention.
\label{page:two-sided}

\subsection{Continued fractions} One writes
$Y_0=\CF{a_0\\a_1\\a_2\\\ldots}$, where formally
\begin{equation}
\label{eq:definition}
\CF{a_0\\a_1\\a_2\\\ldots\\a_h}=a_0+1/\,\CF{a_1\\a_2\\\ldots\\a_{h-1}}
\quad\text{and}\quad \CF{\phantom{a}}=\infty\,.
\end{equation}
It follows, again by induction on $h$, that the definition 
$$\begin{pmatrix}a_0&1\\0&1\end{pmatrix}
\begin{pmatrix}a_1&1\\0&1\end{pmatrix}\cdots
\begin{pmatrix}a_h&1\\0&1\end{pmatrix}=:
\begin{pmatrix}x_h&x_{h-1}\\y_h&y_{h-1}\end{pmatrix}
$$
entails $\CF{a_0\\a_1\\a_2\\\ldots\\a_h}=x_h/y_h$. This provides a
correspondence between the \emph{convergents} $x_h/y_h$ and
certain products of $2\times2$ matrices (more precisely, between the
sequences $(x_h)$, $(y_h)$ of \emph{continuants} and those matrices).
It is a useful exercise to notice that
$Y_0=\CF{a_0\\a_1\\\ldots\\a_h\\Y_{h+1}}$ implies that
$$
Y_{h+1}=-(y_{h-1}Y-x_{h-1})/(y_{h}Y-x_{h})
$$
and that this immediately gives
\begin{equation}\label{eq:distance}
Y_1Y_2\cdots Y_{h+1}=(-1)^{h}(x_h-y_h Y)^{-1}.
\end{equation}
The quantity $-\deg(x_h-y_h Y)=\deg y_{h+1}$ is a weighted sum
giving a measure of the `distance' traversed by the
\cfe\ to its $(h+1)$-st complete quotient. Taking norms yields
\begin{equation}
\label{eq:norm}
(x_h-y_h Y)(x_h-y_h\conj Y)=(-1)^{h+1}Q_{h+1}\,.
\end{equation}

\setcounter{footnote}{0}

\subsection{Conjugation, symmetry, and periodicity}\label{ss:csp}
Each
\pq\ $a_h$ is the \poly\ part of its
corresponding \cq\ $Y_h$. Note, however, that the assertions
above are independent of that conventional selection rule.

One readily shows that $Y_0$ \emph{reduced}, to wit $ \deg Y_0\gt0$ and $\deg\ov
Y_0\lt0$, implies that each complete quotient
$Y_h$ is reduced.  Indeed, it also follows that
$\deg B_h \gt0$, while plainly $\deg\conj
B_h\lt0$ since
$-\conj B_h$ is a remainder; so the $B_h$ too are reduced. In
particular $a_h$, the \poly\ part of $Y_h$, is also the
\poly\ part of $B_h$.

Plainly, at least the first two leading terms of each \poly\ $P_h$
must coincide with the leading terms of $Y-T$. It also follows that
the
\poly s
$P_h$ and
$Q_h$ satisfy the bounds
\begin{equation}\label{eq:bounds} \deg P_h=g+1
\quad\text{and}\quad
\deg Q_h\le g\,.
\end{equation}
Thus, if the base field $\F$ is finite the box
principle entails the
\cfe\ of $Y_0$ is periodic.  If $\F$ is infinite, periodicity is
just happenstance. 

Suppose, however, that the function field
$\F(X,Y)$ is exceptional in that $Y_0$, say, does have a periodic \cfe.
If the \cfe\ of $Y_0$ is periodic
then, by conjugation, also the \ex\ of $B_0$ is periodic.
But conjugation reverses the order of the lines
comprising a
\cf\ tableau. Hence the conjugate of any  preperiod is a
`postperiod', an absurd notion. It follows that,  if periodic, the
two conjugate
\ex s are purely periodic.

Denote by $A$ the \poly\ part of $Y$, and recall that $Y+\conj
Y=T$. It happens that line~$0$ of the \cfe\ of $Y+A-T$ is  
\begin{equation} \label{eq:line0}
Y+A-T=2A-T-(\conj Y+A-T)
\end{equation}
and is symmetric. In general, if  the
\ex\ of
$Y_0$ has a symmetry, and if the \cfe\ is
periodic, its period must have a second symmetry\footnote{The case
of period length $1$ is an exception unless we count its one line as having two
symmetries; alternatively unless we deem it to have period $r=2$.}. So if $Y$ is
exceptional in having a periodic \cfe\ then its period is of
lenght $2s$ and has an additional symmetry of the first kind
$P_s=P_{s+1}$, or its period is of length $2s+1$ and also has a 
symmetry of the second kind, $Q_s=Q_{s+1}$. Conversely, this is
the point, if the \ex\ of $Y$ has a second symmetry then it must be
periodic as just described.
\label{page:periodic}

\subsection{Units}\label{ss:units} It is easy to apply the Dirichlet box principle
to prove that an order
$\Q[\omega]$ of a quadratic number field $\Q(\omega)$ contains
nontrivial units. 
Indeed, by that principle there are infinitely
many pairs of integers $(p,q)$ so that $|q \omega-p|\lt1/q$, whence
$|p^2-(\omega+\conj \omega)pq+ \omega\conj \omega q^2|\lt
(\omega-\conj \omega)+1$. It follows, again by the box principle, 
that there is an integer
$l$ with $0\lt |l|\lt (\omega-\conj \omega)+1$ so that the equation
$p^2-(\omega+\conj \omega)pq+ \omega\conj \omega q^2=l$ has
infinitely many pairs $(p,q)$ and $(p',q')$ of solutions with
$p\=p'$ and $q\=q'\pmod l$. For each such distinct pair,
$xl=pp'-\omega\conj \omega qq'\!$,
$yl=pq'-p'q+(\omega+\conj \omega)qq'$, yields
$(x-\omega y)(x-\conj \omega y)=1$.

In the function field case, we cannot apply the the box principle
for a second time if the base field $\F$ is infinite. So the
existence of a nontrivial unit $x(X)-y(y)Y(X)$ is exceptional.
This should not be a surprise. By the definition of the notion
`unit', such a unit
$u(X)$ say, has a divisor supported only at infinity. Moreover,
$u$ is a function of the order $\F[X,Y]$, and is say of degree
$m$, so the existence of $u$ implies that the class
containing the divisor at infinity is a torsion divisor on
the Jacobian of the curve~\eqref{eq:eq}. The existence of such a
torsion divisor is of course exceptional.

 Suppose now that the function field $\F(X,Y)$ does contain a
nontrivial unit
$u$, say of norm
$-\kappa$ and degree~$m$. Then $\deg(yY-x)=-m\lt-\deg y$, so $x/y$
is a convergent of $Y$ and so some
$Q$ is
$\pm
\kappa$, say
$Q_r=\kappa$ with $r$ odd. That is, line~$r$ of the \cfe\ of $Y+A-T$
is
\begin{equation}\label{eq:liner}
Y_r:=(Y+A-T)/\kappa=2A/\kappa-(\conj
Y+A-T)/\kappa\,;\tag*{line~r:}\end{equation}  here we have used
the fact that $(Y+P_r)/\kappa$ is reduced to deduce that necessarily
$P_r=P_{r+1}=A-T$. 

By conjugation of the $(r+1)$-line tableau commencing with
\eqref{eq:line0} we see that  
\begin{equation}\label{eq:line2r}
Y_{2r}:=Y+A-T=2A-T-(\conj Y+A-T)\,,\tag*{line~2r:}\end{equation}
so that in any case if $Y+A-T$ has a quasi-periodic \cfe\ then it
is periodic of period twice the quasi-period. This result of
Tom~Berry \cite{Be} applies to arbitrary quadratic irrationals
with \poly\ trace. Other elements $(Y+P)/Q$ of
$\F(X,Y)$, with $Q$ dividing the norm   $(Y+P)(\conj Y+P)$, may
be honest-to-goodness quasi-periodic, that is, not also periodic. 

\label{page:quasiperiodic}

Further, if $\kappa\ne-1$ then $r$ \emph{must} be odd. To see that,
notice the identity
\begin{equation}\label{eq:Schmidt}
B\CF{Ca_0\\Ba_1\\Ca_2\\Ba_3\\\ldots}=C\CF{Ba_0\\Ca_1\\Ba_2\\Ca_3\\\ldots},
\end{equation}
reminding one how to multiply a \cfe\ by some
quantity; this cute formulation of the multiplication rule is due
to Wolfgang Schmidt \cite{Schm}.  The `twisted symmetry' occasioned
by division by $\kappa$, equivalent to the existence of a
non-trivial quasi-period, is noted by Christian Friesen \cite{Fr}. 

In summary, if the \cfe\ of $Y$ is quasi-periodic it is periodic,
and the \ex\ has the symmetries of  the more familiar number field
case, as well as  twisted symmetries occasioned by a nontrivial
$\kappa$.

One shows readily that if $x/y=\CF{A\\a_1\\\ldots\\a_{r-1}}$
and $x-Y y$ is a unit of the domain $\F[X, Y]$ then, with $a_{r-1}=\kappa a_1$, 
$a_{r-2}= a_2/\kappa$, $a_{r-3}=\kappa a_3$, $\ldots\,$,
$$\CF{\overline{2A-T\\a_1\\\ldots\\a_{r-1}\\(2A-T)/\kappa\\a_{r-1}\\\ldots\\a_1}}$$
is the quadratic irrational Laurent series $Y+A-T$. Alternatively,
given the \ex\ of $Y+A-T$, and noting that therefore
$\deg Q_r=0$, the fact that the said \ex\ of $x/y$ yields a unit
follows directly from
\eqref{eq:norm}.

\section{Comments}\label{s:comments}

\subsection{} According to Gauss (\emph{Disquisitiones Arithmetic\ae},
Art. 76) \ldots\ \emph{veritates ex notionibus potius quam ex hauriri
debebant}\footnote{[mathematical] truths flow from notions
rather than from notations.}. Nonetheless, one should not underrate the importance
of notation; good notation can decrease the viscosity of
the flow to truth. From the foregoing it seems clear that, given $Y^2=A^2+4v(X-w)$,
one should study the
\cfe\ of $Z=\frac12(Y+A)$, as is done in~\cite{AR}. Moreover, it is a mistake
to be frustrated by minimal models $V^2+UV-vV=U^3-fU+vwU$.

Specifically, we understand that $V^2-8vV=U^3-(4f-1)U^2+8v(2w-1)U$ yields
$Y^2=(X^2+4f-1)^2+4\cdot8v\bigl(X-(2w-1)\bigr)$ by way of $2U=X^2+Y+(4f-1)$ and 
$(V-8v)=XU$. Now $X\gets 2X+1$, $Y\gets 4Y$ means that, instead, we obtain 
$Y^2=(X^2+X+f)^2+4v\bigl(X-(w-1)\bigr)$. This derives from $V^2+UV-vV=U^3-fU+vwU$
by taking $2U=X^2+X+Y+f$ and $V-v=XU$.

\subsection{} The discussion above may have some interest for
its own sake, but my primary purpose is to test ideas for generalisation to higher
genus~$g$. An important difficulty when $g\gt1$ is that \pq s may be of
degree greater than one without that entailing periodicity, whence my 
eccentric aside at page~\pageref{a:aside}. Happily, the generalisation to
translating by a point $(w_0, e_0-e_1)$ on the quartic model effected above also is
a simplification in that one surely may always choose a translating divisor so as
to avoid meeting singular steps in the \cfe. 
In that context one
finds that the sequence $(\ldots\,, 2, 1,1,1,1,1, 1,
2,3,4,8, 17,50, 
\ldots\,)$ satisfying the recursion
$T_{h-3}T_{h+3}=T_{h-2}T_{h+2}+T_{h}^2$  arises from adding multiples of the
divisor at infinity on the Jacobian of  the curve
$Y^2=(X^3-4X+1)^2+4(X-2)$ of genus~$2$ to the divisor 
$[(\varphi,0),(\conj\varphi,0)]$; $\varphi$ is the golden ratio.

\bibliographystyle{amsalpha}

\label{page:lastpage}
\end{document}